\documentclass[a4paper,twoside,11pt,reqno]{amsart}

\usepackage[T1]{fontenc}
\usepackage{courier}
\usepackage[scaled=0.92]{helvet}

\usepackage[english]{babel}
\usepackage{amsthm,amsfonts,amssymb,amsmath}
\usepackage{hyperref}
\hypersetup{
colorlinks=true,
citecolor=red,
linkbordercolor={1 1 1},
citebordercolor={1 1 1},
}

\usepackage{mathrsfs}
\usepackage{draftwatermark}
\SetWatermarkLightness{1} 
\SetWatermarkAngle{45}
\SetWatermarkText{DRAFT}
\usepackage{xcolor} 

\theoremstyle{plain}
\newtheorem{satz}{Theorem}[section]
\newtheorem{prop}[satz]{Proposition}

\newtheorem{lem}[satz]{Lemma}

\theoremstyle{definition}

\newtheorem{hyp}[satz]{Hypothesis}

\newcommand{\la}{\langle}
\newcommand{\re}{\rangle}
\newcommand{\mx}{\mbox}
\newcommand{\rw}{\rightarrow}
\newcommand{\de}{\displaystyle}
\newcommand{\ml}{\mathcal}

\newcommand{\pl}{\partial}

\newcommand{\x}{\times}

\newcommand{\ue}[1]{\underline{#1}}
\newcommand{\uue}[1]{\underline{\underline{#1}}}
\newcommand{\uuue}[1]{\underline{\underline{\underline{#1}}}}
\newcommand{\beq}[1]{\begin{equation} \label{#1}}
\newcommand{\eeq}{\end{equation}}
\newcommand{\beqar}{\[ \begin{array}{rcl}}
\newcommand{\eeqar}{\end{array} \]}
\newcommand{\beqdesar}{\[ \begin{array}{rcll}} 
\newcommand{\eeqdesar}{\end{array} \]}

\providecommand{\ep}{\varepsilon}
\providecommand{\ph}{\varphi}

\providecommand{\RR}{\mathbb{R}}
\providecommand{\CC}{\mathbb{C}}

\providecommand{\ZZ}{\mathbb{Z}}
\providecommand{\NN}{\mathbb{N}}

\providecommand{\TT}{\mathbb{T}}

\newcommand{\snorm}[2]{\left| #1\right|_{#2}}
\newcommand{\norm}[2]{\left \lVert#1 \right\rVert_{#2}}

\usepackage[utf8]{inputenc}

\hypersetup{pdfauthor={Fortunati, Wiggins},%
            pdftitle={Kolmogorov theorem for Aperiodic Poisson Systems},%
            pdfsubject={Time-dependent Poisson systems},%
            pdfkeywords={Poisson systems, Aperiodic time dependence}%
}

\DeclareMathOperator{\id}{Id}

\usepackage{pstricks,pst-node,pst-tree}
\psset{npos=.6,nrot=:U,labelsep=2pt,tnpos=a,treesep=0.5cm}

\makeatletter
\g@addto@macro{\endabstract}{\@setabstract}
\newcommand{\authorfootnotes}{\renewcommand\thefootnote{\@fnsymbol\c@footnote}}%
\makeatother

\begin{document}
\title[A Kolmogorov theorem for nearly-integrable Poisson systems...]{A Kolmogorov theorem for nearly-integrable Poisson systems with asymptotically decaying time-dependent perturbation.}
\author{Alessandro Fortunati}
\thanks{This research was supported by ONR Grant No.~N00014-01-1-0769 and MINECO: ICMAT Severo Ochoa project SEV-2011-0087.}
\address{School of Mathematics, University of Bristol, Bristol BS8 1TW, United Kingdom}
\email{alessandro.fortunati@bristol.ac.uk}

\author{Stephen Wiggins}
\email{s.wiggins@bristol.ac.uk}

\maketitle

\begin{abstract}
The aim of this paper is to prove the Kolmogorov theorem of persistence of Diophantine flows for nearly-integrable Poisson systems associated to a real analytic Hamiltonian with aperiodic time dependence, provided that the perturbation is asymptotically vanishing. The paper is an extension of an analogous result by the same authors for canonical Hamiltonian systems; the flexibility of the Lie series method developed by A. Giorgilli et al., is profitably used in the present generalisation. 
\end{abstract}
\section{Preliminaries}
Since the theory of Hamiltonian system has been extended to more general spaces than the standard symplectic manifold (see e.g. \cite{laucavan}), the interest for stability theory for these ``generalized'' (or ``non-canonical'') systems has rapidly risen. Among them, examples of systems defined on \emph{Poisson manifolds}, also known as \emph{Poisson systems}, naturally appear in applications. Fundamental stability results such as the Nekhoroshev and Kolmogorov theorems can be obtained for such class of systems, see \cite{liyi} and \cite{conghong}.\\        
The main goal of this paper is to extend the proof of the Kolmogorov theorem for aperiodically time-dependent canonical systems considered in \cite{forwig} to a class of Poisson system admitting a family of invariant tori, as considered in \cite{liyi}. The (necessary) assumption of slow dacay of the perturbation is still assumed in this case while restriction to the quadratic dependence in the actions of \cite{forwig} is removed. 
\\
The proof, closely carried out along the lines of \cite{forwig} (which is the extension of \cite{gior} to the time-dependent case), is based on the tools of the Lie series method used in \cite{gior} and \cite{nuovocimento}, after a straightforward extension of them to the Poisson case.  
\subsection{Poisson systems: a short summary}
Let $Z:=Y \times \TT^{n}$ where $Y \subset \RR^m$, with $m,n \in \NN$. Denoting $\ue{z} \equiv (\ue{x},\ue{y}) \in Z$, 
we consider the following system of ODEs
\beq{eq:poisson}
\ue{\dot{z}}=\uue{\mathcal{B}}(\ue{z}) \ml{H}_{\ue{z}}(\ue{z}) \mx{,}
\eeq
where\footnote{we denote with $\ml{M}_{N,M}(U;\RR)$ the space of $N \times M$ real (complex, if $\CC$) valued matrix defined on $U$.} $\uue{\mathcal{B}}(z)=\{b_{kl}(z)\} \in \mathcal{M}_{(n+m),(n+m)}(Z;\RR)$ and $\ml{H}:Z \rw \RR$ (Hamiltonian) are given. By defining the \emph{brackets} of two function $F,G:Z \rw \RR$ as
\beq{eq:brackets}
\{F,G\}^*:=(F_{\ue{z}})^T \uue{\mathcal{B}} \cdot G_{\ue{z}}  \mx{,}
\eeq
system (\ref{eq:poisson}) reads as $\dot{\ue{z}}=\{\ue{z},\ml{H}\}^*$. \\
The brackets (\ref{eq:brackets}) are said to be \emph{Poisson brackets} if $\uue{\mathcal{B}}$ is such that $\{\cdot,\cdot\}^*$ is a \emph{skew-symmetric}, \emph{bilinear} form satisfying: $\{F G, H\}^*=F\{G,H\}^*+G \{F,H\}^*$ and $\{\{F,G\},H\}^*+\{\{G,H\},F\}^*+\{\{H,F\},G\}^*=0$ for all $F,G,H$, i.e. the Leibnitz and Jacobi identities, respectively. See \cite{marrat} for a comprehensive treatment. Correspondingly, (\ref{eq:poisson}) is called \emph{Poisson system} on $Z$.\\
As in the canonical case, if $\ml{H}=\ml{H}(\ue{z},t)$ with $t \in [0,+\infty) $, the obtained time-dependent Poisson system can be interpreted as an autonomous system in an extended space. More precisely, by setting $\xi:=t$, considering the new Hamiltonian $H(\ue{z},\eta,\xi):=\eta+\ml{H}(\ue{z},\xi)$ and the matrix
\[
\uue{B}:=
\left(
\begin{array}{cc}
\uue{\ml{B}} & \uue{0} \\ 
\uue{0} & \uue{J}
\end{array}
\right), 
\qquad
\uue{J}:=
\left(
\begin{array}{cc}
0 & -1 \\ 
1 & 0
\end{array}
\right),
\]
one obtains the following equivalent system
\beq{eq:poissontime}
\left(
\begin{array}{c}
\dot{\ue{z}} \\
\dot{\eta}\\
\dot{\xi}
\end{array}
\right)
=
\uue{B} H_{(\ue{z},\eta,\xi)} \mx{.}
\eeq
where $(\ue{z},\eta,\xi) \in Z \x \RR \x \RR^+ =: \ml{D}$, if the evolution of the auxiliary variable $\eta$ is disregarded.\\ 
It is straightforward to check that the brackets associated to $\uue{B}$, i.e. defined by
\beq{eq:bracketstime}
\{F,G\}:=(F_{(\ue{z},\eta,\xi)})^T \uue{B} \cdot G_{(\ue{z},\eta,\xi)} \mx{,}
\eeq
for all $F,G:\ml{D} \rw \RR$, satisfy the properties described above, i.e. system (\ref{eq:poissontime}) is a Poisson system on $\ml{D}$.

\section{Framework and main result}
\subsection{Nearly-integrable Poisson systems and their invariant tori}
From now on we shall consider time-dependent and nearly-integrable Poisson systems, i.e. associated to Hamiltonians of the form
\beq{eq:hamnotempo}
H(\ue{y},\ue{x},\eta,\xi):=h(\ue{y})+\eta+ \ep f(\ue{y},\ue{x},\eta,\xi) \mx{.}
\eeq
As discussed in \cite{liyi}, the skew-symmetry condition and the necessity to obtain an invariant torus once a particular value $\ue{y}^* \in Y$ has been chosen, implies that the matrix $\uue{\ml{B}}$ has the following particular structure
\[
\uue{\ml{B}}=\uue{\ml{B}}(\ue{y})=
\left(
\begin{array}{cc}
\uue{0} & \uue{\ml{B}}_{12}\\
-\uue{\ml{B}}_{12}^T & \uue{\ml{B}}_{22} 
\end{array}
\right),
\] 
with $\uue{\ml{B}}_{12} \in \ml{M}_{m,n}(Y; \RR)$ and $\uue{\ml{B}}_{22} \in \ml{M}_{n,n}(Y; \RR)$ is skew-symmetric. In such a way the unperturbed vector field has equations
\[
\ue{\dot{y}}=0,\qquad \ue{\dot{x}}=\ue{\omega}, \qquad \dot{\eta}=0, \qquad \dot{\xi}=1 \mx{,}
\]
where the vector
\beq{eq:omega}
\ue{\omega}:=\uue{\ml{B}}^0 \ue{\tilde{\omega}}; \qquad 
\uue{\mathcal{B}}^0 := -\uue{\ml{B}}_{12}^T(\ue{y}^*), \qquad \ue{\tilde{\omega}}:=h_{\ue{y}}(\ue{y}^*)  \mx{,}
\eeq
is the frequency of the flow on the chosen invariant torus.\\
In the typical scenario of a Kolmogorov-type result, our aim is to show that the motion with frequency $\ue{\omega}$ persists in the perturbed system, provided that $\ep$ is ``sufficiently'' small, under suitable hypotheses on the Hamiltonian (\ref{eq:hamnotempo}).\\
After a (formal) Taylor expansion of the Hamiltonian (\ref{eq:hamnotempo}) around $\ue{y}^*$, a trivial rescaling of the expansion point in the origin and a redefinition of $\ue{y}$ and of $H$, the same Hamiltonian reads, up to an additive constant as
\beq{eq:hamstart}
H(\ue{y},\ue{x},\eta,\xi)=\ue{\tilde{\omega}} \cdot \ue{y}+ \eta + \frac{1}{2} \uue{\ml{C}} \ue{y} \cdot \ue{y}+ \ep f(\ue{y},\ue{x},\eta,\xi)  + \ml{R}(\ue{y}) \mx{,}
\eeq    
where $\uue{\ml{C}}:= h_{\ue{y} \ue{y}}(\ue{y}^*)$ and $\ml{R}=O(|y|^3)$.
\subsection{Setting and main statement}
In order to use the tools of Complex Analysis, let us consider the \emph{complexification} of the space $\{\ue{0}\} \times \TT^n \x \RR \x \RR^+$, defined by $\ml{D}_{\rho,\sigma;\zeta}:=\Delta_{\rho} \times \TT_{\sigma}^n \times \ml{S}_{\rho} \times \ml{R}_{\zeta}$ where 
$$
\begin{array}{rclrcl}
\Delta_{\rho}&:=&\{\ue{y} \in \CC^m:|\ue{y}| \leq \rho\},
& \qquad 
\TT_{\sigma}^n&:=&\{\ue{x} \in \CC^n: |\Im \ue{x}| \leq \sigma\},\\
\ml{S}_{\rho}&:=&\{\eta \in \CC: |\Im \eta| \leq \rho\}, & \qquad 
\ml{R}_{\zeta}&:=&\{\xi =:\xi_R+i\xi_I \in \CC:\xi_R \geq -\zeta ; \,|\xi_I| \leq \zeta
\} \mx{,}
\end{array}
$$
Analogously to \cite{gior}, the space of scalar valued functions $g=g(\ue{y},\ue{x},\xi)$ defined on $\ml{D}_{\rho,\sigma;\zeta}$ is endowed with the usual \emph{supremum} and \emph{Fourier} norms
\[
\snorm{g}{[\rho,\sigma;\zeta]}:=\sup_{(\ue{y},\ue{x}) \in \ml{D}_{\rho,\sigma;\zeta}} |g(\ue{y},\ue{x},\xi)| \mx{,}
\qquad 
\norm{g}{[\rho,\sigma;\zeta]}:=\sum_{\ue{k} \in \ZZ^n} \snorm{g_{\ue{k}}(\ue{y},\xi)}{[\rho; \zeta]} 
e^{|\ue{k}|\sigma} \mx{,}
\]
where $g_k(\ue{y},\xi)$ are the coefficients of the Fourier expansion $g=\sum_{k \in \ZZ^n} g_k(\ue{y},\xi) e^{i \ue{k} \cdot \ue{x}}$ and $|\ue{k}|:=|k_1|+\ldots+|k_n|$ for all $\ue{k} \in \ZZ^n$. In the case of vector-valued functions $\ue{w}:\ml{D} \rw \CC^l $ we denote $\norm{\ue{w}}{[\rho,\sigma;\zeta]}:=\sum_{j=1}^l \norm{w_j}{[\rho,\sigma;\zeta]}$. Given a matrix $\uue{M} \in \ml{M}_{n,m}(\Delta_{\rho},\CC)$, the following norm will be finally considered: $\norm{\uue{M}}{\rho}:=n m \max_{ij} (\sup_{\ue{y} \in \Delta_{\rho}}|m_{ij}(y)|)$ while we shall set simply $\norm{\uue{M}}{}:=n m \max_{ij} |m_{ij}|$ if $\uue{M}$ does not depend on $\ue{y}$.\\
The function $h$ and the matrix $\uue{\ml{B}}(y)$ will be supposed analytic and bounded on $\ml{D}_{\rho,\sigma;\zeta}$, i.e. there exist two constants $M_{h},M_{\mathcal{\uue{B}}}>0$ such that $\norm{h}{\rho} \leq M_{h}$ and $\norm{\uue{\mathcal{B}}}{\rho} \leq M_{\mathcal{B}}$. In particular the expansion leading to (\ref{eq:hamstart}) is well defined.\\
Analogously to \cite{forwig} we shall assume the following hypothesis
\begin{hyp}\label{hyp}
\begin{itemize}
\item There exists $\upsilon \in (0,1)$ such that, for all $ \ue{v} \in \CC^m$ 
\beq{eq:hypongamma}
 |\uue{\ml{C}} \ue{v}| \leq {\upsilon}^{-1}|\ue{v}| \mx{.}
\eeq
\item (Diophantine condition): $\ue{y}^*$ is such that $\ue{\omega}$ is a $\gamma-\tau$ Diophantine vector, i.e. there exist $\gamma$ and $\tau>n-1$ such that $|\la \ue{\omega}, \ue{k} \re | \geq \gamma |\ue{k}|^{-\tau}$, for all $\ue{k} \in \ZZ^{n}\setminus \{\ue{0}\}$.
\item (Slow decay): the perturbation $f$ is holomorphic on $\ml{D}_{\rho,\sigma;\zeta}$, satisfying, in addition,  
\beq{eq:slowdecay}
\norm{f(\ue{y},\ue{x},\xi)}{[\rho,\sigma/2;\zeta]} \leq M_f e^{-a |\xi|} \mx{,}
\eeq
for some $M_f>0$ and $a \in (0,1)$. 
\end{itemize}
\end{hyp}
The choice $a<1$ is not technical but simply related to the possibility to obtain simpler estimates in the follow. Nevertheless, it allows us to exploit the slow decaying feature of the perturbation.\\
In the described framework, we are able to prove the following
\begin{satz}[Aperiodic-Poisson Kolmogorov]\label{thm}
Consider (\ref{eq:hamstart}) under the Hypothesis \ref{hyp}.\\ Then, for all $a \in (0,1)$ there exists $\ep_a>0$ such that for all $\ep \in (0,\ep_a]$, it is possible to find an analytic, $\ep-$close to the identity, Poisson change of variables $(\ue{y},\ue{x},\eta,\xi) = \ml{P}(\ue{y}^{(\infty)},\ue{x}^{(\infty)},\eta^{(\infty)},\xi)$, $\ml{P}:\ml{D}^{*} \subset \ml{D} \rw \ml{D}$ 
casting Hamiltonian (\ref{eq:hamstart}) into the \emph{Kolmogorov normal form} 
\beq{eq:kolnormal}
H^{(\infty)}(\ue{y}^{(\infty)},\ue{x}^{(\infty)},\eta^{(\infty)},\xi)=\ue{\tilde{\omega}} \cdot  \ue{y}^{(\infty)} + \eta^{(\infty)}+
\mathcal{R}_{\infty}(\ue{y}^{(\infty)},\ue{x}^{(\infty)},\xi;\ep) \mx{,}
\eeq
where $\ml{R}_{\infty}$ is at least quadratic in $\ue{y}$.
\end{satz}

\section{Formal scheme}
As usual, the construction of the Kolmogorov normal form is based on an iterative algorithm of changes of variables. We shall use the \emph{Lie method} by considering the \emph{Lie series operator} associated to the \emph{generating function} $\phi$, formally defined as
\[
\exp(\ml{L}_{\phi})=\id + \sum_{s \geq 1} \ml{L}_{\phi}^s
\] 
where $\ml{L}_{\phi}\cdot:=\{\phi,\cdot\}$ the latter being defined in (\ref{eq:bracketstime}). By construction, the above defined operator, is a Poisson change, as a time-one evolution of the Poisson system associated to the ``Hamiltonian'' $\phi$.\\
The hard-core aspect of the normalization algorithm consists in the following
\begin{lem}\label{lem:formalscheme}
Let us suppose that, for some $j \in \NN$, Hamiltonian (\ref{eq:hamstart}) has the form
\beq{eq:hamricorsiva}
H^{(j)}=\eta + \tilde{h}^{(j)} \mx{,}
\eeq
where $\tilde{h}^{(j)}=h^{(j)}+ g^{(j)}$ with
\beq{eq:hgj}
h^{(j)}:=\de \tilde{\ue{\omega}}\cdot\ue{y}+\frac{1}{2}\uue{C}^{(j)}\ue{y}\cdot\ue{y}+R^{(j)}, \qquad
g^{(j)}:= \de A^{(j)}+\ue{B}^{(j)} \cdot \ue{y} \mx{,}
\eeq
with $A^{(j)}, \ue{B}^{(j)}, \uue{C}^{(j)}$ and $R^{(j)}$ depending on $\ue{x},\xi$. Moreover $R^{(j)}=O(|\ue{y}|^3)$.\\
Then it is possible to find $\chi^{(j)}$ such that $\ml{P}_j H^{(j)}$ 
with 
\beq{eq:pj}
\ml{P}_j:=\exp(\ml{L}_{\chi^{(j)}}) \mx{,}
\eeq
is of the form (\ref{eq:hamricorsiva}) for suitable $A^{(j+1)},\ue{B}^{(j+1)},\uue{C}^{(j+1)}$ and $R^{(j+1)}$.
\end{lem}
The effect of this scheme is to remove the presence of the ``unwanted'' terms collected in $g^{(j)}$ on a certain ``level'' of magnitude\footnote{The use of a book-keeping parameter (see e.g. \cite{efthy}) could be very useful in order to to recognize the perturbative feature of the scheme. Set e.g. $H^{(j)}=\eta + h^{(j)}+ \lambda_j g^{(j)}$ (then, by writing $H^{(0)}$ is easy to see that $\lambda_0=O(\ep)$) and repeat the computation below considering the operator $\exp(\ml{L}_{\lambda_j\chi^{(j)}})$. This will show that the terms removed by the homological equation lie on the level $\lambda_j$ and that one can set $\lambda_{j+1}=O(\lambda_j^2)$, exploiting the well known quadratic feature of the Kolomogorov scheme. Unfortunately, the use of this parameter is not particularly effective in the quantitative part and it will be avoided.}. This cancellation is effected via the time-dependent homological equation (\ref{eq:homological}), already introduced in \cite{forwig}. The transformation determined in this way produces further terms of this type i.e. $A^{(j+1)}$ and $\ue{B}^{(j+1)}$, but their size is ``smaller'' than the same terms labelled with $j$. In this way, their contribution is (formally) removed once one sets 
\beq{eq:compmap}
\ml{P}:=\lim_{j \rw \infty} \ml{P}_{j} \circ \ml{P}_{j-1} \circ \ldots \circ \ml{P}_1 \mx{.}
\eeq
This well established heuristic approach will be made rigorous in the quantitative part.
\proof
We consider the action of the operator $\exp(\ml{L}_{\chi^{(j)}})$ on $H_j$, obtaining
\beqar
\hat{H}^{(j)}:=\exp(\ml{L}_{\chi^{(j)}})H^{(j)}&=& \eta+ h^{(j)}+ g^{(j)} +  \chi_{\xi}^{(j)}+ \{ \chi^{(j)},h^{(j)} \} \\
&+& \{ \chi^{(j)},g^{(j)} \}+\de  \sum_{s \geq 2} \frac{1}{s!}
\left[ 
\ml{L}_{\chi^{(j)}}^s \eta+ \ml{L}_{\chi^{(j)}}^s \tilde{h}^{(j)} \right]\mx{.}
\eeqar
Let us suppose that it is possible to determine $\chi^{(j)}$ such that
\beq{eq:homological}
\chi_{\xi}^{(j)}+g^{(j)}+\{ \chi^{(j)},h^{(j)} \} = Q^{(j)}(\ue{y},\ue{x},\xi;)=O(|\ue{y}|^2) \mx{.}
\eeq
In such case the Hamiltonian takes the form
\[
\hat{H}^{(j)}= \de \eta + h^{(j)} +Q^{(j)}+ \hat{R}^{(j)}; \qquad
\hat{R}^{(j)}:=\{ \chi^{(j)},g^{(j)} \}+\sum_{s \geq 2} \frac{1}{s!}
\ml{L}_{\chi^{(j)}}^s \left(  \eta+  \tilde{h}^{(j)} \right) \mx{.}
\]
Note that $h^{(j)}+\lambda_j Q^{(j)}$ is at least quadratic in $\ue{y}$ while $\hat{R}^{(j)}$ contains also terms independent of $\ue{y}$ and linear in $\ue{y}$. Hence, it is possible to set 
\begin{subequations}
\begin{align}
A^{(j+1)}:=& \de  \hat{R}^{(j)}(\ue{0}) \label{eq:ajp}\\
\ue{B}^{(j+1)}:=& \de   \hat{R}_{\ue{y}}^{(j)}(\ue{0}) \label{eq:bjp}\\
\uue{C}^{(j+1)}:=& \de  \hat{H}_{\ue{y}\ue{y}} ^{(j)} (\ue{0}) \equiv \uue{C}^{j}+ (\exp(\ml{L}_{\chi^{(j)}})\tilde{h}^{(j)}-\tilde{h}^{(j)})_{\ue{y} \ue{y}}(\ue{0}) \label{eq:cjp}\mx{.}
\end{align}
\end{subequations}
By using (\ref{eq:hgj}) one defines $h^{(j+1)}$ and $g^{(j+1)}$. Then $H^{(j+1)}$ by (\ref{eq:hamricorsiva}). The residual higher order terms of the Taylor expansion are stored in $R^{(j+1)}$.
\endproof
\subsection{Solution of the Homological equation}
Our aim is now to determine a solution of the equation (\ref{eq:homological}). Recalling (\ref{eq:bracketstime}), equation (\ref{eq:homological}) takes the form 
\beq{eq:newhomological}
 \chi_{\xi}^{(j)}+A^{(j)}+\ue{B}^{(j)} \cdot \ue{y} 
-\chi_{\ue{x}}^{(j)} \uue{\mathcal{B}}_{12}^T \cdot (\ue{\omega}+\uue{C}^{(j)} \ue{y}) = O(|\ue{y}|^2) \mx{.}
\eeq
The necessity to solve the previous equation up to first order in $\ue{y}$, leads to the possibility to restrict ourselves to linear expansions of $\uue{\mathcal{B}}$ and to the well known class of linear generating functions, as suggested by Kolmogorov 
\[
\uue{\mathcal{B}}_{12}^T  = -\uue{\mathcal{B}}^0-\uuue{\mathcal{B}}^1\ue{y}+O(|\ue{y}|^2), \qquad \chi^{(j)}=S^{(j)}+\ue{T}^{(j)}\cdot \ue{y} \mx{,}
\] 
having recalled (\ref{eq:omega}) and set $\uuue{\mathcal{B}}^1:=-(\uue{B}_{12}^{T})_{\ue{y}}(\ue{y}^*)$. Plugging these expansions into (\ref{eq:newhomological}), the comparison of the power of $\ue{y}$ yields, up to $O(|\ue{y}|^2)$, the following system
\beq{eq:system}
\left\{
\begin{array}{rcl}
S_{\xi}^{(j)}+S_{\ue{\omega}}^{(j)} +A^{(j)} & = & 0\\
\ue{T}_{\xi}^{(j)}+ \ue{T}_{\ue{\omega}}^{(j)} +S_{\ue{x}}^{(j)} \uue{E}^{(j)} + \ue{B}^{(j)} & = & 0
\end{array}
\right. \mx{.}
\eeq
Having denoted $\uue{E}^{(j)}:=\uue{\mathcal{B}}^0 \uue{C}^{(j)} +  \uuue{N}^1 \tilde{\ue{\omega}} $ and $\pl_{\ue{\omega}}:=\pl_{\ue{x}} \cdot \ue{\omega}$. \\
Once $S^{(j)}$ has been determined by solving the first equation, each component of the second equation has exactly the same form of the first one, with the corresponding component of $S_{\ue{x}}^{(j)} \uue{E}^{(j)} + \ue{B}^{(j)}$ in place of $A^{(j)}$. This completes the formal resolvability of the iterative step.\\
We stress that equations of the system (\ref{eq:system}), both of the form
\beq{eq:homologicaltime}
\ph_{\xi} + \ph_{\ue{\omega}} = \psi \mx{,}
\eeq
possess the same structure of those found in the canonical case, discussed in \cite{forwig}.

\section{Quantitative estimates on the formal scheme}
\subsection{Technical tools}
The following two statements are excerpted from \cite{forwig} and reported below for the reader's convenience, as they play a key role in the quantitative part of the proof. The first one concerns a bound for the composition of an arbitrary number of Lie operators. The second one provides a result of existence and analyticity for the solution of the time-dependent homological equation (\ref{eq:homologicaltime}). 

\begin{prop}\label{prop:chipsi}
Let $d_1,d_2 \in [0,1/2]$ and $\chi$ and $\psi$ be two functions on $\ml{D}_{\rho,\sigma;\zeta}$ such that $\norm{\chi}{[(1-d_1)(\rho,\sigma);\zeta]}$ and $\norm{\psi}{[(1-d_2)(\rho,\sigma);\zeta]}$ are bounded for all $\xi \in \ml{R}_{\zeta}$.\\
Then for all $\tilde{d} \in (0,1-\hat{d})$ where $\hat{d}:=\max\{d_1,d_2\}$ and for all $s \geq 1$ one has the following estimate  
\beq{eq:iterativo}
\norm{\ml{L}_{\chi}^s \psi}{[(1-\tilde{d}-\hat{d})(\rho,\sigma);\zeta]} \leq 
\frac{s!}{e^2} \left( \frac{4e^2 \Gamma_{\rho,\sigma}}{\tilde{d}^2}\right)^s
\norm{\chi}{[(1-d_1)(\rho,\sigma);\zeta]}^s 
\norm{\psi}{[(1-d_2)(\rho,\sigma);\zeta]} \mx{.}
\eeq
where $\Gamma_{\rho,\sigma}:=[e^2 G_{11}  \sigma^2+ 2 e G_{12} \rho\sigma  
+ G_{22} \rho^2]( e \rho \sigma )^{-2}$ with $G_{ij}:=\norm{\uue{\mathcal{B}}_{ij}}{\rho}$.
\end{prop}
\proof Straightforward generalisation of \cite[Pag. 77]{giorgilli02}.
\endproof
It is immediate to see that, under the same assumptions, the Lie operator $\exp(\ml{L}_{\chi})$ converges if the following condition is satisfied
\beq{eq:convergence}
\mathfrak{L}:= \frac{4e^2 \Gamma_{\rho,\sigma}}{\tilde{d}^2} \norm{\chi}{[(1-d_1)(\rho,\sigma);\zeta]} \leq \frac{1}{2} \mx{.}
\eeq
\begin{prop}\label{prop:small} Let $\delta \in [0,1)$ and suppose that there exists a constant $K>0$ such that 
\beq{eq:expdec}
\norm{\psi}{[(1-\delta)\hat{\sigma};\hat{\zeta}]} \leq K e^{-a |\xi|} \mx{,}
\eeq
for all $\xi \in \ml{R}_{\hat{\zeta}}$ and for some $0<\hat{\sigma}<\sigma$, $0<\hat{\zeta} \leq \zeta$. Note that $a$ has been defined in (\ref{eq:slowdecay}).\\
Then for all $d \in (0, 1-\delta)$ and for all $\hat{\zeta}$ such that
\beq{eq:sceltazeta}
4 |\ue{\omega}| \hat{\zeta} \leq d \hat{\sigma} \mx{,}
\eeq
the solution of (\ref{eq:homologicaltime}) exists and satisfies
\beq{eq:stimehom}
\norm{\ph}{[(1-\delta-d)\hat{\sigma};\hat{\zeta}]}  \leq  
\frac{K \Theta_1}{a(d \hat{\sigma})^{2\tau}} e^{-a |\xi|}, \qquad
\norm{\ph_{x_l} }{[(1-\delta-d)\hat{\sigma};\hat{\zeta}]}  
\leq \frac{K \Theta_2}{a (d \hat{\sigma})^{2\tau+1}} 
e^{-a |\xi|} \mx{,}
\eeq
for all $l=1,\ldots,n$ and for some constants $\Theta_2>\Theta_1>0$. 
\end{prop}
\proof Given in \cite{forwig}. \endproof

\subsection{Iterative lemma}
Let us define the following vector of parameters $\ue{u}_j:=(d_j,\epsilon_j,\zeta_j,m_j,\rho_j,\sigma_j) \in [0,1)^6$ for all $j \geq 0$. Consider, in addition, $\ue{u}_*:=(0,0,0,m_*,\rho_*,\sigma_*)$ for some $m_*,\rho_*,\sigma_* \in (0,1)$. The vectors $\ue{u}_*,\ue{u}_0$ will be determined later. We shall denote $\ml{D}^{(j)}:=\ml{D}_{\rho_j,\sigma_j;\zeta_j}$ and $\ml{D}^*:=\ml{D}_{\rho_*,\sigma_*;\zeta_*}$.
\begin{lem}\label{lem:iterative}
Under the same hypothesis of Lemma \ref{lem:formalscheme}, suppose, in addition, the existence of $\ue{u}_j$ with $\ue{u}_j>\ue{u}_*$ (i.e. component-wise), such that the following conditions hold true
\begin{enumerate}
\item 
\beq{eq:iterativeitemone}
\max\left\{\norm{A^{(j)}}{[\sigma_j;\zeta_j]},\norm{\ue{B}^{(j)}}{[\sigma_j;\zeta_j]}\right\} \leq \epsilon_j e^{-a|\xi|}
\mx{,}
\eeq
\item for all functions $\ue{w}=\ue{w}(q,\xi) : \ml{D}^{(j)} \rw \CC^m$ 
\beq{eq:iterativeitemthree}
\norm{\uue{C}^{(j)}(q,\xi) \ue{w}(q,\xi)}{[\sigma_j;\zeta_j]}  \leq \upsilon_j^{-1} \norm{\ue{w}(q,\xi)}{[\sigma_j;\zeta_j]} \mx{,}
\eeq
\item $d_j \leq 1/6$ and $\zeta_j$ is set as
\beq{eq:zetaj}
4 |\ue{\omega}| \zeta_j=d_j \sigma_j \mx{,}
\eeq
\item there exists a constant $M_{\tilde{h}^{(j)}}>0$ such that 
\beq{eq:limh}
\norm{\tilde{h}^{(j)}}{[\rho_j,\sigma_j;\zeta_j]} \leq M_{\tilde{h}^{(j)}} \mx{.}
\eeq
\end{enumerate}
Then it is possible to determine a constant $D>0$ such that: if 
\beq{eq:piccolaunmezzo}
\epsilon_j \frac{D}{a^4 \upsilon_j^2 d_j^{8(\tau+1)}} \leq \frac{1}{2} \mx{,}
\eeq
it is possible to choose $\ue{u}_{j+1} < \ue{u}_j$  under the constraint (\ref{eq:zetaj})\footnote{I.e. satisfying $4 |\ue{\omega}| \zeta_{j+1}=d_{j+1} \sigma_{j+1}$.}, for which (\ref{eq:iterativeitemone}), (\ref{eq:iterativeitemthree}), (\ref{eq:limh}) are satisfied by $A^{(j+1)},\ue{B}^{(j+1)}$, $\uue{C}^{(j+1)}$ and $\tilde{h}^{(j+1)}$ given by (\ref{eq:ajp}), (\ref{eq:bjp}), (\ref{eq:cjp}) and (\ref{eq:hgj}) respectively, with $M_{\tilde{h}^{(j+1)}} = M_{\tilde{h}^{(j)}}$.
\end{lem}
The proof of this result is organized in the following three steps. In order to avoid a cumbersome notation, the index $j$ will be dropped from all the objects depending on it, and reintroduced only for objects at the $j+1-$th stage. 

\subsubsection{Estimates on the generating function}
By (\ref{eq:iterativeitemone}) and Prop. \ref{prop:small} (set $\delta=d/2$) we get
\beq{eq:sj}
\norm{S}{[(1-d/2)\sigma;\zeta]} \leq  \epsilon  \frac{M_0}{a d^{2\tau}}e^{-a |\xi|},\qquad
\norm{S_{\ue{x}}}{[(1-d/2)\sigma;\zeta]} \leq \epsilon \frac{M_1}{a d^{2\tau+1}}e^{-a |\xi|} \mx{,}
\eeq
where $
M_0:=\Theta_1(2/\sigma_*)^{2\tau}$ and $M_1:=n \Theta_2(2/\sigma_*)^{2\tau+1}$. In this way, recalling the definition of $\uue{E}$, the symmetry of $\uue{C}$, using the second equation of (\ref{eq:sj}) and finally (\ref{eq:iterativeitemthree}), we obtain
\beq{eq:boh}
\norm{S_{\ue{x}} \uue{E} +\ue{B}}{[(1-d/2)\sigma;\zeta]}
\leq 
 \epsilon
\frac{M_2}{a \upsilon d^{2\tau+1}} e^{-a|\xi|}
\eeq
where $M_2:=1+ M_1 \left( 
\norm{\uue{\mathcal{B}}^0 }{}
+
\norm{\uuue{\mathcal{B}}^1 \ue{\omega}}{}
\right) 
$. By (\ref{eq:boh}) and Prop. \ref{prop:small} for the second equation of (\ref{eq:system}) hold 
\beq{eq:bohdue}
\norm{\ue{T}}{[(1-d)\sigma;\zeta]} \leq \epsilon \frac{M_3}{a^2 \upsilon d^{4 \tau+1}} e^{-a|\xi|},\qquad
\norm{\ue{T}_{\ue{x}}}{[(1-d)\sigma;\zeta]} \leq 
\epsilon \frac{M_4}{a^2 \upsilon d^{4 \tau+2}} e^{-a|\xi|}
\eeq
with $M_3:=m M_2 \Theta_1(2/\sigma_*)^{2\tau}$ and $M_4:=m n M_2 \Theta_2(2/\sigma_*)^{2\tau+1}$.\\
Recalling the definition of $\chi$, (\ref{eq:bohdue}) imply the following estimates
\beq{eq:stimechi}
\norm{\chi}{[(1-d)(\rho,\sigma);\zeta]} \leq \epsilon \frac{M_5}{a^2 \upsilon d^{4 \tau+2}} e^{-a|\xi|}, \qquad
\norm{\chi_{\ue{x}}}{[(1-d)(\rho,\sigma);\zeta]}  \leq \epsilon \frac{M_6}{a^2 \upsilon d^{4 \tau+2}} e^{-a|\xi|} \mx{,} 
\eeq
with $M_5:=M_0+M_3$ and $M_6:=M_1+M_4$.
By a Cauchy estimate we immediately get
\beq{eq:chixi}
\norm{\chi_{\xi}}{[(1-d)(\rho,\sigma;\zeta)]} \leq
\frac{1}{d \zeta} \norm{\chi}{[(1-d)(\rho,\sigma);\zeta]}
\leq \epsilon \frac{M_5}{a^2 v d^{4 \tau+3} \zeta} e^{-a|\xi|} \mx{.}
\eeq
\subsubsection{Estimates on the transformed Hamiltonian}
Our aim is to give an estimate for $\hat{R}$. First of all note that by (\ref{eq:iterativo}) written for $s-1$ where\footnote{we set $d_2=d$ and $d_1=d$, then $\hat{d}=d$ and $\tilde{d}=d$ as $\tilde{d}$ needs to be smaller than $1-2d$. This is true as we shall suppose $d \leq 1/6$.} $\psi \equiv \ml{L}_{\chi} \eta = \chi_{\xi} $ (and $\chi$ as itself), one gets, using (\ref{eq:stimechi}) and (\ref{eq:chixi})
\beq{eq:chieta}
\norm{\ml{L}_{\chi}^s \eta}{[(1-2 d)(\rho,\sigma;\zeta)]} 
\leq  \frac{s!}{4 s e^4 \Gamma_{\rho,\sigma}\zeta} \mathfrak{L}^s,
\qquad \mathfrak{L}:=\epsilon  \frac{4 e^2 \Gamma_{\rho ,\sigma} M_5}{a^2 \upsilon d^{4 \tau+3}} e^{-a|\xi|}  \mx{.}
\eeq
Hence, the Lie operator is uniformly convergent provided that
\beq{eq:smallone}
\epsilon  \frac{8 e^2 M_{\mathcal{B}} M_5}{a^2 \upsilon d^{4 \tau+3} (\rho_* \sigma_*)^2} \leq \frac{1}{2} \mx{,}
\eeq
as $\Gamma_{\rho ,\sigma} \leq 2 M_{\mathcal{B}}(\rho_* \sigma_*)^{-2}$, see Prop. \ref{prop:chipsi}. Then $\mathfrak{L} \leq 1/2$. The obtained bounds yield
\beq{eq:lieeta}
\sum_{s \geq 2} \frac{1}{s!} \norm{\ml{L}_{\chi}^s \eta}{[(1-2 d)(\rho,\sigma;\zeta)]} \leq \frac{\mathfrak{L}^2}{2 e^2 \Gamma_{\rho,\sigma} \zeta} \mx{.}
\eeq
An estimate of the third term appearing in $\hat{R}$ can be obtained in the same way. More precisely, bound (\ref{eq:iterativo}) whith $\psi \equiv 
\tilde{h}$ yields 
$ \norm{\ml{L}_{\chi}^s \tilde{h}}{[(1-2 d)(\rho,\sigma;\zeta)]} \leq 
M_{\tilde{h}} e^{-2} s! \mathfrak{L}^s $, yielding
\beq{eq:third}
\sum_{s \geq 2} \frac{1}{s!} \norm{\ml{L}_{\chi}^s \tilde{h}}{[(1-2 d)(\rho,\sigma;\zeta)]} \leq  \frac{2 M_{\tilde{h}}}{e^2}  \mathfrak{L}^2 \mx{.}
\eeq
Finally, the first term of $\hat{R}$ follows easily from (\ref{eq:iterativo}) for $s=1$. Indeed, as $\norm{g}{[\rho,\sigma;\zeta]} \leq 2 \epsilon e^{-a|\xi|} \mx{,}
$ by (\ref{eq:iterativeitemone}), we get $\norm{\ml{L}_{\chi}^s g}{[(1-2 d)(\rho,\sigma;\zeta)]} \leq \epsilon 2 e^{-2} \mathfrak{L}e^{-a|\xi|}$. The latter, (\ref{eq:third}) and (\ref{eq:lieeta}) imply
\[
\norm{\hat{R}}{[(1-2 d)(\rho,\sigma;\zeta)]} \leq \epsilon^2 \frac{M_7}{a^4 \upsilon^2 d^{8 \tau+6}\zeta} e^{-2 a |\xi|} \mx{,}
\]
with $M_7:=16 M_{\mathcal{B}} M_5 (1 + 8 e^2 M_{\mathcal{B}} M_{\tilde{h}} M_5 + e^2 M_5)(\rho_* \sigma_*)^{-4} $, allowing us to obtain the following bounds 
\begin{subequations}
\begin{align}
\norm{A^{(j+1)}}{[(1-3 d)(\rho,\sigma;\zeta)]} & \leq \de 
\epsilon^2 \frac{M_7}{a^4 \upsilon^2 d^{8 \tau+6}\zeta} e^{-2 a |\xi|} \mx{,} \label{eq:ajpuno}
\\ 
\norm{\ue{B}^{(j+1)}}{[(1-3 d)(\rho,\sigma;\zeta)]} & \leq  \de 
\frac{1}{d \rho}\norm{\hat{R}}{[(1-2 d)(\rho,\sigma;\zeta)]} \leq 
\epsilon^2 \frac{M_7}{a^4 \upsilon^2 d^{8 \tau+7}\zeta \rho_*} e^{-2 a |\xi|} \mx{.} \label{eq:bjpuno}
\end{align}
\end{subequations}
The final step is the estimate of $\uue{C}^{(j+1)}$. Taking into account  (\ref{eq:cjp}) and the bound before (\ref{eq:third}) we compute
\[
\norm{\exp(\ml{L}_{\chi}) \tilde{h}-\tilde{h}}{[(1-2 d)(\rho,\sigma;\zeta)]}  \leq  \de \sum_{s \geq 1} \frac{1}{s!} \norm{\ml{L}_{\chi}^s \tilde{h}}{[(1-2 d)(\rho,\sigma;\zeta)]}  \leq \epsilon \frac{8 \Gamma_{\rho,\sigma} M_{\tilde{h}} M_5}{a^2 v d^{4 \tau+3}\zeta} e^{-a|\xi|} \mx{.}
\]
Now denoting $\uue{C}':=\uue{C}^{(j+1)}$ and $M_8:=32 m M_{\mathcal{B}} M_{\tilde{h}} M_5 (\rho_*^2 \sigma_*)^{-2}$ we obtain
\beq{eq:ckl}
\norm{C_{kl}'-C_{kl}}{[(1-3 d)(\rho,\sigma;\zeta)]}
\leq \frac{2}{(d\rho)^2} \norm{\exp(\ml{L}_{\chi}) \tilde{h}-\tilde{h}}{[(1-2 d)(\rho,\sigma;\zeta)]} \leq \epsilon \frac{ M_8}{a^2 m v d^{4 \tau+5}\zeta} e^{-a|\xi|}
\eeq
In conclusion we have, for all $\ue{w}=\ue{w}(\ue{x},\xi) \in \CC^m$ (write $\uue{C}'=\uue{C}+(\uue{C}'-\uue{C})$ and use (\ref{eq:ckl})),
\beq{eq:cjpu}
\norm{\uue{C}'  \ue{w}}{[(1-3 d)(\rho,\sigma;\zeta)]} \leq 
\left[\upsilon - \epsilon \frac{ M_8}{a^2 v d^{4 \tau+5}\zeta} e^{-a|\xi|} \right]^{-1} \norm{\ue{w}}{[(1-3 d)(\sigma;\zeta)]}
=: \de (\upsilon')^{-1} \norm{\ue{w}}{[(1-3 d)(\sigma;\zeta)]} \mx{.}
\eeq
The previous definition of $\upsilon'$ makes sense provided that the quantity between the square brackets is positive, which is a property that can be obtained by requiring 
\beq{eq:smallv}
 \epsilon \frac{ M_8}{a^2 v^2 d^{4 \tau+5}\zeta}  \leq \frac{1}{2} \mx{.}
\eeq
\subsubsection{Parameters of the iteration}
It is now sufficient to take into account conditions (\ref{eq:smallone}), (\ref{eq:smallv}) and those obtained comparing (\ref{eq:ajpuno}) and (\ref{eq:bjpuno}) with (\ref{eq:iterativeitemone}), to see immediately that these condition hold by (\ref{eq:piccolaunmezzo}) for a suitable\footnote{One can set e.g. $D:=32 e^2|\ue{\omega}| M_{\mathcal{B}} (\rho_* \sigma_*)^{-2} \max\{M_6,M_7,M_8\}$ (as $M_7>M_5$ and by (\ref{eq:zetaj})). See also Sec \ref{sec:change}. } $D$ and setting
\beq{eq:quad}
\epsilon_{j+1} := \frac{D}{a^4 \upsilon_j^2 d_j^{8(\tau+1)}}  \epsilon_j^2 \mx{,}
\eeq
which is the well known quadratic iteration, that is able to compensate the effect of the small divisors contained in $d_j^{-8(\tau+1)}$. The monotonicity of $\epsilon_j$ is a direct consequence as $\epsilon_0$ will be chosen smaller than one.\\
Taking into account the domain restrictions appearing in (\ref{eq:ajpuno}), (\ref{eq:bjpuno}) and (\ref{eq:cjpu}), the iterative step is complete once 
the following conditions are set
\beq{eq:paramiter}
(\rho_{j+1},\sigma_{j+1})=(1-3 d_j)(\rho_j,\sigma_j), \qquad \upsilon_{j+1}=(1-d_j^{4 \tau + 3})\upsilon_j \mx{,}
\eeq
where the latter follows directly from (\ref{eq:quad}). As for $\zeta_{j+1}$, it is easy to see that condition (\ref{eq:zetaj}) for $j+1$ is stronger than the restriction required by the above mentioned bounds\footnote{It is sufficient to show that $\hat{\zeta}_{j+1}:=(4|\ue{\omega}|)^{-1} d_{j+1} \sigma_{j+1}$ (i.e. (\ref{eq:zetaj}) for $j+1$) is smaller than $\tilde{\zeta}_{j+1}:=(1-3d_j)\zeta_j$. Using (\ref{eq:zetaj}) in the latter and using the just obtained value for $\sigma_{j+1}$ one has that $\hat{\zeta}_{j+1}< \tilde{\zeta}_{j+1}$ provided that $d_{j+1} < d_j$, a property that will be made true by construction.}.\\ 
Analogously to \cite{nuovocimento}, the property (\ref{eq:limh}), follows from the fact that $\exp(\ml{L}_{\chi^{(j)}})$ maps points $\ue{z} \in \mathcal{D}^{(j+1)} \subset \mathcal{D}^{(j)}$ to $\exp(\ml{L}_{\chi^{(j)}}) \ue{z}=:\ue{z}' \in \mathcal{D}^{(j)}$ and from the well known relation $\exp(\ml{L}_{\chi^{(j)}}) 
\hat{h}^{(j)}(\ue{z})=\hat{h}^{(j)}(\ue{z}')$. This completes the proof of Lemma \ref{lem:iterative}. 
\subsection{Bounds on the Poisson change of coordinates}\label{sec:change}
As in \cite{forwig}, the aim is now to give the estimates for the change of variables $\ml{P}_{j}$ for all $j \geq 0$. These will be used later to ensure that the image of points in $\ml{D}^*$ do not ``escape'' from $\ml{D}^{(0)}$ under the action of the composition of Poisson flows $\ml{P}_{j}$ once the limit (\ref{eq:compmap}) is taken. 
\begin{prop} In the hypotheses of Lemma \ref{lem:iterative} the Poisson change of variables $\ue{z}^{(j)}=\ml{P}_j (\ue{z}^{(j+1)})$ is an $\epsilon_0-$close to the identity change of variables, satisfying $|\xi^{(j+1)}-\xi^{(j)}|=0$ i.e. $\xi^{(j)} \equiv \xi $ for all $j \geq 1$ and 
\beq{eq:bounds}
\max\{|\ue{y}^{(j+1)}-\ue{y}^{(j)}|,|\eta^{(j+1)}-\eta^{(j)}|\}  
\leq  d_j \rho_j e^{-a |\xi|}, \qquad 
|\ue{x}^{(j+1)}-\ue{x}^{(j)}|  \leq  d_j \sigma_j e^{-a |\xi|} \mx{.}
\eeq
\end{prop}
\proof Firstly note that $\ml{L}_{\chi} \ue{y} = \de - \chi_{\ue{x}} \uue{\ml{B}}_{12}^T$, $\ml{L}_{\chi} \ue{x} = \chi_{\ue{x}} \uue{\ml{B}}_{22} + \ue{T} \uue{\ml{B}}_{12}$ and $\ml{L}_{\chi} \xi = 0$, in particular $\xi^{j+1}=\xi^{j} \equiv \xi$ (the Poisson transformation does not act on time). By using the bounds of the previous section it is possible to bound $\norm{\ml{L}_{\chi} \ue{y}}{[(1-d)(\rho,\sigma;\zeta)]}$ and $
\norm{\ml{L}_{\chi} \ue{x}}{[(1-d)(\rho,\sigma;\zeta)]}$. A repeated use of (\ref{eq:iterativo}) written for $s-1$, with $\ml{L}_{\chi} \ue{y}$ and $\ml{L}_{\chi} \ue{x}$ in place of $\psi$, respectively, yields
\[
\norm{\ml{L}_{\chi}^s \ue{y}}{[(1-2 d)(\rho,\sigma;\zeta)]} 
\leq \de \frac{s! M_{\mathcal{B}} M_6 }{4 s e^4 \Gamma_{\rho,\sigma} M_5} \mathfrak{L}^s, \qquad 
\norm{\ml{L}_{\chi}^s \ue{x}}{[(1-2 d)(\rho,\sigma;\zeta)]}  \leq \frac{s! M_{\mathcal{B}} (M_3+M_6)}{4 s e^4 \Gamma_{\rho,\sigma} M_5} \mathfrak{L}^s \mx{,} 
\]
while the bound for $\ml{L}_{\chi}^s \eta$ is already known from (\ref{eq:chieta}). 
 The $\epsilon_0-$ closeness of the change of variables is evident by the previous bounds and by (\ref{eq:quad}).\\
By construction $\ue{z}^{(j)}=\exp(\ml{L}_{\chi^{(j)}})\ue{z}^{(j+1)}$, hence $|\ue{z}^{(j+1)}-\ue{z}^{(j)}| \leq \sum_{j \geq 1} (s!)^{-1}\norm{\ml{L}_{\chi^{(j)}}^s \ue{z}^{(j+1)}}{[(1-2 d)(\rho,\sigma;\zeta)]}$, in this way the use of (\ref{eq:piccolaunmezzo}) and (\ref{eq:zetaj}) give the desired estimates.
\endproof
\section{Convergence of the scheme and conclusion of the proof}
This part is completely similar to \cite{forwig} with the exception of some details, which are specified in the brief sketch reported below for completeness.\\
Set $\epsilon_j=\epsilon_0(j+1)^{-16 (\tau+1)}$ for all $j \geq 1$ and 
\beq{eq:dj}
d_j= \left( \frac{D \epsilon_0}{a^4 \upsilon_j^2} \right)^{\frac{1}{8(\tau+1)}} \frac{(j+2)^2}{(j+1)^4} \mx{.}
\eeq
Under the condition 
\beq{eq:epz}
\epsilon_0 \leq \frac{a^4}{D 12^{8(\tau+1)}} \mx{,}
\eeq
holds $d_j < d_0 \leq 1/6$ as required by Lemma {\ref{lem:iterative}}, item (3). Furthermore, due to (\ref{eq:dj}), the sequences $\rho_j,\sigma_j,\upsilon_j$ as determined by (\ref{eq:paramiter}), converge to some $\rho_*,\sigma_*,\upsilon_*>0$, respectively. More precisely one finds $(\rho_*,\sigma_*,\upsilon_*)=(\rho_0/4,\sigma_0/4,\upsilon_0/2)$, with $\rho_0,\sigma_0$ and $\upsilon_0$ to be determined. By (\ref{eq:dj}), (\ref{eq:zetaj}) and the definition of $\epsilon_j$ above, it follows that $\lim_{j \rw \infty} \ue{u}_j=\ue{u}_*$.\\
Consider the Taylor expansion of Hamiltonian (\ref{eq:hamstart}) and set $A^{(0)}(\ue{x},\xi):=\ep f(\ue{x},\ue{0},\xi)$, $\ue{B}^{(0)}(\ue{x},\xi):=\ep [f_{\ue{y}}(\ue{x},\ue{0},\xi)]$ and $\uue{C}^{(0)}:=\uue{\mathcal{C}}+ \ep [f_{\ue{y} \ue{y}}(\ue{x},\ue{0},\xi)]$, finally storing in $R^{(0)}$ the higher order terms. The expanded Hamiltonian is immediately in the form (\ref{eq:hamricorsiva}) and it is possible to set $H^{(0)}:=H$.\\
By assumption (\ref{eq:slowdecay}) and a Cauchy estimate one has $\norm{f_{\ue{y}} (\ue{x},\ue{y},\xi)}{[\rho_0,\sigma_0;\zeta_0]} \leq m M_f \rho_0^{-1} \exp(-a|\xi|)  $ setting $(\rho_0,\sigma_0)=(\rho/2,\sigma/2)$. In this way (\ref{eq:iterativeitemone}) holds true by choosing $\epsilon_0:=m \ep M_f/\rho_0$. Setting $\upsilon_0=\upsilon/2$ 
condition (\ref{eq:iterativeitemthree}) is satisfied for sufficiently small $\ep$ (see \cite{forwig} for a quantitative estimate). The described conditions on $\ep$, together with (\ref{eq:epz}), determine the threshold for $\ep_a$ mentioned in theorem \ref{thm}. More precisely, by (\ref{eq:epz}), it is of the form $\ep_a \leq M_a O(a^4)$, where $M_a$ is a (``very small'') constant. Note that the property (\ref{eq:limh}) for $H^{(0)}$ is a direct consequence of the regularity assumptions on the initial Hamiltonian. For instance, one can set $M_{\hat{h}^{(0)}}:=M_{h}+\ep M_f$. The value of $\zeta_0$ is determined by (\ref{eq:zetaj}) by setting $d_0=1/6$. The choice of $\ue{u}_0$ as a function of the initial domain of analyticity is now complete. \\
Finally, bounds (\ref{eq:bounds}) and (\ref{eq:dj}) under the condition (\ref{eq:epz}), ensure that e.g. $|\ue{y}^{(\infty)}-\ue{y}| \leq \sum_{j \geq 0} |\ue{y}^{(j+1)}-\ue{y}^{(j)}| \leq d \rho/6$ and similarly for the other coordinates, hence points starting in $\ml{D}^{*}$ are mapped within $\ml{D}^{(0)}$. The degeneration of the radius  $\zeta$ is not relevant as the transformation is trivial in time. The proof of the analyticity of $\ml{P}$, as defined in (\ref{eq:compmap}), follows by the Weierstra{\ss} theorem, see e.g. \cite[Pag. 168]{gior}.

\bibliographystyle{alpha}
\bibliography{noncanKAM.bib}

\end{document}